\documentclass[11pt]{article}
\usepackage{amsmath}
\usepackage{amssymb}
\newtheorem{theorem}{\hspace*{\parindent}Theorem}
\newtheorem{lemma}{\hspace*{\parindent}Lemma}

\newtheorem{conjecture}{\hspace*{\parindent}Conjecture}
\title{Positivity of Toeplitz determinants formed by rising factorial series and properties of related polynomials}
\author{D.\,Karp\footnote{Institute of Applied Mathematics, Vladivostok, Russia,
e-mail:\,\emph{dimkrp@gmail.com}}}
\begin{document}
\maketitle

\begin{center}
\parbox{12cm}{
\small\textbf{Abstract.} In this note we prove positivity of
Maclaurin coefficients of polynomials written in terms of rising
factorials and arbitrary log-concave sequences.  These polynomials
arise naturally when studying log-concavity of rising factorial
series.  We propose several conjectures concerning  zeros and
coefficients of a generalized form of those polynomials.  We also
consider polynomials whose generating functions are higher order
Toeplitz determinants formed by rising factorial series. We make
three conjectures about these polynomials.  All proposed
conjectures are supported by numerical evidence.}
\end{center}

\bigskip

Keywords: \emph{Log-concavity, P\'{o}lya frequency sequences,
Toeplitz determinant, stability, hyperbolicity, rising factorial,
hypergeometric functions}

\bigskip

MSC2010:  26C10, 05A20

\bigskip

\renewcommand{\thetheorem}{\Alph{theorem}}
\setcounter{theorem}{0}

\paragraph{1. Introduction.}  The confluent hypergeometric
function is defined by the series
\begin{equation}\label{eq:1F1defined}
{_{1}F_1}(a;c;z)
:=\sum\limits_{n=0}^{\infty}\frac{(a)_n}{(c)_n}\frac{z^n}{n!},
\end{equation}
where $(a)_0=1$, $(a)_n=a(a+1)\cdots(a+n-1)=\Gamma(a+n)/\Gamma(a)$
is rising factorial or Pochhammer symbol.  It was proved by
Barnard, Gordy and Richards in \cite{BGR} that the function
$$
z\to\begin{vmatrix}
{_{1}F_1}(a;c;z) & {_{1}F_1}(a+1;c;z) \\
{_{1}F_1}(a-1;c;z) & {_{1}F_1}(a;c;z)
\end{vmatrix}
$$
has positive Maclaurin coefficients if $a>0$, $c>-1$ ($c\neq{0}$).
This has been extended by Karp and Sitnik in \cite{KS} to the
determinant ($\alpha,\beta>0$)
\begin{equation}\label{eq:f-defined}
z\to\begin{vmatrix}
f(x+\alpha;z) & f(x+\alpha+\beta;z) \\
f(x;z) & f(x+\beta;z)
\end{vmatrix},
~\text{where}~f(x;z):=\sum\limits_{n=0}^{\infty}f_k(x)_k\frac{z^k}{k!}
\end{equation}
and $\{f_k\}_{k=0}^{\infty}$ is any non-negative log-concave
sequence without internal zeros, i.e. $f_k^2\geq{f_{k-1}f_{k+1}}$,
$k=1,2,\ldots,$ and if $f_N=0$ for some $N>0$ then $f_k=0$ for all
$k\geq{N}$. Since
$$
\begin{vmatrix}
f(x+\alpha;z) & f(x+\alpha+\beta;z) \\
f(x;z) & f(x+\beta;z)
\end{vmatrix}=\sum\limits_{n=2}^{\infty}Q^{\alpha,\beta}_{n}(x)\frac{z^n}{n!},
$$
where
\begin{equation}\label{eq:Qn}
Q^{\alpha,\beta}_{n}(x):=\sum\limits_{k=0}^{n}f_{k}f_{n-k}\binom{n}{k}\left[(x+\alpha)_k(x+\beta)_{n-k}-(x+\alpha+\beta)_k(x)_{n-k}\right],
\end{equation}
Theorem~1 from \cite{KS} can be restated as follows:
\begin{theorem}\label{th:KS}
Suppose $\{f_{k}\}_{k=0}^{n}$ is a non-negative log-concave
sequence without internal zeros, $\alpha,\beta>0$, $n\geq{2}$.
Then $Q^{\alpha,\beta}_n(x)\geq{0}$ for all $x\geq{0}$. The
inequality is strict unless $f_k=q^k$, $k=0,1,\ldots,n$ for some
$q>0$.
\end{theorem}
Note that this theorem does not cover the above result from
\cite{BGR} completely since Theorem~\ref{th:KS} requires $x$ to be
non-negative while the result in \cite{BGR} is valid for
$x=a-1>-1$. On several occasions (see, for instance,
\cite{AIM,KCMFT} I proposed the following two conjectures:
\begin{conjecture}\label{cj:positivity}
If $f_k^2>f_{k-1}f_{k+1}$, $k=1,2,\ldots,n-1$, $n\geq{3}$, then
$Q^{\alpha,\beta}_n(x)$ has positive coefficients at $x^j$,
$j=0,1,\ldots,n-2$
\end{conjecture}

Recall that a polynomial is called Hurwitz stable if all its zeros
have negative real part.  See details and extensions in \cite{T}.
\begin{conjecture}\label{cj:hurwitz}
If $f_k^2>f_{k-1}f_{k+1}$, $k=1,2,\ldots,n-1$, $n\geq{3}$, then
$Q^{\alpha,\beta}_n(x)$ is Hurwitz stable.
\end{conjecture}
For polynomials with real coefficients stability implies
positivity of coefficients (this result is usually attributed to
A.\,Stodola (1893)) so that Conjecture~\ref{cj:positivity} is true
if Conjecture~\ref{cj:hurwitz} holds.

Let me also propose a third conjecture that has not been presented
elsewhere. It requires the notion of P\'{o}lya frequency sequence
defined  formally in section~4 below.  Briefly,
$\{f_k\}_{k=0}^{n}$ is $PF_{\infty}$ if all minors of the infinite
matrix (\ref{eq:Polya}) are non-negative.
\begin{conjecture}\label{cj:zeros}
If $\{f_k\}_{k=0}^{n}$ is $PF_{\infty}$, $n\geq{3}$, then all
zeros of $Q^{1,1}_n(x-1)$ are real and negative.
\end{conjecture}

Let me remark that Conjecture~\ref{cj:zeros} fails for
$Q^{\alpha,\beta}_n(x)$ with arbitrary $\alpha,\beta>0$ and so
does it for $Q^{1,1}_n(x-1)$ when $\{f_k\}_{k=0}^{n}$ is only
log-concave ($PF_{\infty}$ is much stronger requirement than
log-concavity, see details in section~4). I have explicit (but a
bit cumbersome) counterexamples that demonstrate these claims. All
three conjectures are supported by massive numerical evidence.

In a relatively recent work \cite{IsLaf} Ismail and Laforgia and,
more recently, Baricz and Ismail \cite{BarIsm} proved absolute or
complete monotonicity of numerous Hankel determinants formed by
special functions which possess the integral representation
$$
f_n=\int\limits_{\alpha}^{\beta}[\phi(t)]^nd\mu(t),
$$
where both the function $\phi$ and the measure $\mu$ may depend on
parameters.  When the size of the determinant is equal to 2 their
results reduce to positivity and integral representations for
$f_{n}f_{n+2}-f_{n+1}^2$. The positivity of this expression is
discrete log-convexity of (or reverse Tur\'{a}n type inequality
for) $f_n$.  Unfortunately, the technique used in these papers
does not extend to log-concavity  (discrete or not) as far as I
can see, although some discrete log-concavity results are proved
in \cite{BarIsm} employing a different method.

The  purpose of this note is twofold.  First, we prove the
positivity of the coefficients of $Q_{n}^{1,1}(x-1)$ settling a
particular case of Conjecture~\ref{cj:positivity}. This furnishes
a far-reaching extension of the result of \cite{BGR} and partially
of \cite{KS}. Second, we consider a higher order Toeplitz
determinant whose entries are functions defined in
(\ref{eq:f-defined}). We give power series expansion of such
determinant in powers of $z$ with coefficients being polynomials
in $x$.  We make several conjectures about these polynomials
serving as natural generalizations of
Conjectures~\ref{cj:positivity} and \ref{cj:zeros} for
$Q^{1,1}_n(x-1)$.

\renewcommand{\thetheorem}{\arabic{theorem}}
\setcounter{theorem}{0}

\paragraph{2. Preliminaries.} We will need several lemmas which we
present in this section.  We will always assume that the sequence
$\{f_k\}$ is not a zero sequence.
\begin{lemma}\label{lm:fedor}
Suppose $\{f_k\}_{k=0}^{n}$ has no internal zeros and
$f_k^2\geq{f_{k-1}f_{k+1}}$, $k=1,2,\ldots,n-1$. If the real
sequence $M_0,M_1,\ldots,M_{[n/2]}$ satisfying $M_{[n/2]}>0$ and
$\sum_{k=0}^{[n/2]}M_k\geq{0}$ has one change of sign, then
\begin{equation}\label{eq:keysum}
\sum\limits_{0\leq{k}\leq{n/2}}f_{k}f_{n-k}M_k\geq{0}.
\end{equation}
Equality is only attained if $f_k=\alpha^k$, $\alpha>0$, and
$\sum_{k=0}^{[n/2]}M_k=0$.
\end{lemma}
\textbf{Proof.} Suppose $f_k>0$, $k=s,\ldots,p$, $s\geq{0}$,
$p\leq{n}$. Log-concavity of $\{f_k\}_{k=0}^{n}$ clearly implies
that $\{f_{k}/f_{k-1}\}_{k=s+1}^{p}$ is decreasing, so that for
$s+1\leq{k}\leq{n-k+1}\leq{p+1}$
$$
\frac{f_k}{f_{k-1}}\geq\frac{f_{n-k+1}}{f_{n-k}}~\Leftrightarrow~f_{k}f_{n-k}\geq
f_{k-1}f_{n-k+1}.
$$
Since $k\leq{n-k+1}$ is true for all $k=1,2,\ldots,[n/2]$, the
weights $f_{k}f_{n-k}$ assigned to negative $M_k$s in
(\ref{eq:keysum}) are smaller than those assigned to positive
$M_k$s leading to (\ref{eq:keysum}).  The equality statement is
obvious.~~~$\square$

We will use the formula
\begin{equation}\label{eq:symmgenerating}
\prod\limits_{k=1}^{q}(x+a_k)=\sum\limits_{k=0}^{q}e_{q-k}(a_1,\ldots,a_q)x^{k},
\end{equation}
where $e_m(a_1,\ldots,a_q)$ denotes $m$-th elementary symmetric
polynomial,
$$
e_k(a_1,\ldots,a_q)=\sum\limits_{1\leq{j_1}<{j_2}\cdots<{j_k}\leq{q}}a_{j_1}a_{j_2}\cdots{a_{j_k}}.
$$
The key fact about elementary symmetric polynomials that we will
need requires the notion of majorization \cite[Definition A.2,
formula (12)]{MOA}. It is said that $B=(b_1,\ldots,b_q)$ is weakly
supermajorized by $A=(a_1,\ldots,a_q)$ (symbolized by
$B\prec^W{A}$) if
\begin{equation}\label{eq:amajorb}
\begin{split}
& 0<a_1\leq{a_2}\leq\cdots\leq{a_q},~~
0<b_1\leq{b_2}\leq\cdots\leq{b_q},
\\
&\sum\limits_{i=1}^{k}a_i\leq\sum\limits_{i=1}^{k}b_i~~\text{for}~~k=1,2\ldots,q.
\end{split}
\end{equation}

\begin{lemma}\label{lm:major}
Suppose $B\prec^W{A}$. Then
$$
\frac{e_{k}(a_1,\ldots,a_q)}{e_{k-1}(a_1,\ldots,a_q)}\leq\frac{e_{k}(b_1,\ldots,b_q)}{e_{k-1}(b_1,\ldots,b_q)},~~k=1,2,\ldots,q.
$$
\end{lemma}
\textbf{Proof.}   According to \cite[3.A.8]{MOA} $B\prec^W{A}$
implies that $\phi(A)\leq\phi(B)$ if and only if $\phi(x)$ is
Schur-concave and increasing in each variable. Hence, we should
choose
\[
\phi_k(x_1,\ldots,x_q)=\frac{e_{k}(x_1,\ldots,x_q)}{e_{k-1}(x_1,\ldots,x_q)},~~k=1,2,\ldots,q.
\]
Schur-concavity of these functions has been proved by Schur (1923)
- see \cite[3.F.3]{MOA}.  It is left to show that $\phi_k$ is
increasing in each variable.  Due to symmetry we can take $x_1$ to
be variable thinking of $x_2,\ldots,x_q$ as being fixed.  Using
the definition of elementary symmetric polynomials we see that for
$k\geq{2}$
\[
\phi_k(x_1,\ldots,x_q)=\frac{x_1e_{k-1}(x_2,\ldots,x_q)+e_{k}(x_2,\ldots,x_q)}{x_1e_{k-2}(x_2,\ldots,x_q)+e_{k-1}(x_2,\ldots,x_q)}.
\]
So taking derivative with respect to  $x_1$ we obtain
($e_m=e_m(x_2,\ldots,x_q)$ for brevity):
\[
\frac{\partial\phi_k(x_1,\ldots,x_q)}{\partial{x_1}}=\frac{e_{k-1}(x_1e_{k-2}+e_{k-1})-e_{k-2}(x_1e_{k-1}+e_{k})}
{[x_1e_{k-2}+e_{k-1}]^2}=\frac{e_{k-1}^2-e_{k}e_{k-2}}{[x_1e_{k-2}+e_{k-1}]^2}\geq{0}.
\]
Non-negativity holds by Newton's inequalities.~~$\square$

Next lemma is a part of Theorem~\ref{th:KS}.
\begin{lemma}\label{lm:Chu}
Suppose $f_k=1$ for all $k=0,1,\ldots,n$. Then
$Q_{n}^{\alpha,\beta}(x)\equiv{0}$.
\end{lemma}
\textbf{Proof.}  If $f_k=1$ for all $k=0,1,\ldots,n$, then
$Q_{n}^{\alpha,\beta}(x)/n!$ is $n$-th Maclaurin coefficient of
the function
$$
z\to(1-z)^{-x-\alpha}(1-z)^{-x-\beta}-(1-z)^{-x}(1-z)^{-x-\alpha-\beta}\equiv{0}~~\square
$$

\bigskip

\paragraph{3. Main results.}
Introduce the notation
\begin{equation}\label{eq:Pn}
P_{n}(x)=Q_{n}^{1,1}(x-1)=\sum\limits_{k=0}^{n}f_{k}f_{n-k}\binom{n}{k}\left[(x)_k(x)_{n-k}-(x+1)_k(x-1)_{n-k}\right].
\end{equation}
According to Lemma~\ref{lm:Chu} $P_{n}(x)\equiv{0}$ if $f_k=1$ for
all $k=0,1,\ldots,n$. Our main theorem is as follows.
\begin{theorem}\label{th:size2}
If $f_k^2>f_{k-1}f_{k+1}$ for $k=1,2,\ldots,n-1$, then $P_n(x)$
has degree $n-2$ and positive coefficients.
\end{theorem}
\textbf{Proof.} Denote
$$
\Phi_k(x)=2(x)_{k}(x)_{n-k}-(x-1)_{k}(x+1)_{n-k}-(x-1)_{n-k}(x+1)_{k}~\text{for}~k<n-k
$$
and $\Phi_k(x)=(x)_{k}(x)_{n-k}-(x-1)_{k}(x+1)_{n-k}$ for $k=n-k$
(which only happens for even $n$).  Then
$$
P_n(x)=\sum\limits_{0\leq{k}\leq{n/2}}f_kf_{n-k}\binom{n}{k}\Phi_k(x).
$$
Straightforward computation yields
\begin{equation}\label{eq:Phi0}
\Phi_0(x)=-n(n-1)(x+1)_{n-2},
\end{equation}
\begin{equation}\label{eq:Phik}
\Phi_k(x)=(x)_{k-1}(x+1)_{n-k-2}l_k(x),~~1\leq{k}\leq{n/2},
\end{equation}
where
\begin{equation}\label{eq:lk}
l_k(x)=-A_kx+B_k,
\end{equation}
$$
A_k=n(n-1)-4k(n-k),~~B_k=n(n-1)-2k(n-k),~~1\leq{k}<{n/2},
$$
and
$$
A_{n/2}=-n/2,~~B_{n/2}=n(n-2)/4.
$$
These formulas show that $\Phi_k(x)$ has degree $n-2$ for all
$0\leq{k}\leq{n/2}$ and the free term is only present in
$\Phi_0(x)$, where it equals $-n!$, and in  $\Phi_1(x)$, where it
equals $(n-1)!$.  Hence, the free term in $P_n(x)$ is equal to
$$
-f_0f_{n}\binom{n}{0}n!+f_1f_{n-1}\binom{n}{1}(n-1)!=n!(f_1f_{n-1}-f_0f_{n})>0,
$$
and it remains to prove the theorem for the coefficients at  $x^j$
for $j=1,2,\ldots,n-2$. Since for $n=2$ we only have the free term
we can assume that $n\geq{3}$.

Now if $a_{k,j}$ is the coefficient at $x^j$, $j=1,2,\ldots,n-2$,
in $\Phi_{k}(x)$, $k=0,1,\ldots,[n/2]$, then setting
$M_{k,j}=\binom{n}{k}a_{k,j}$ we have according to
Lemma~\ref{lm:Chu}:
$$
\sum_{0\leq{k}\leq{n/2}}M_{k,j}=0,~~~j=1,2,\ldots,n-2.
$$
Formula (\ref{eq:Phi0}) shows that $a_{0,j}<0$ for all
$j=1,2,\ldots,n-2$. Hence, in order to apply Lemma~\ref{lm:fedor}
we only need to demonstrate that the sequence $a_{k,j}$,
$k=0,1,\ldots,[n/2]$ has precisely one change of sign for each
$j=1,2,\ldots,n-2$. We have
$$
\Phi_{1}(x)=(n-1)(-(n-4)x+n-2)(x+1)_{n-3}.
$$
If $n=3$ then this reduces to $2(x+1)$ and we are done, since the
coefficient at $x$ is positive and $[n/2]=1$, so that
$\Phi_{1}(x)$ is the last term.  If $n=4$ than
$\Phi_{1}(x)=6(x+1)$ and $\Phi_{2}(x)=4x(x+1)$ which again proves
the claim for $n=4$.  Hence, we may assume that $n\geq{5}$.

Formula (\ref{eq:symmgenerating}) and the definition of the
Pochhammer symbol $(x)_m=x(x+1)\cdots(x+m-1)$ lead to
representation
$$
\Phi_k(x)=(x)_{k-1}(x+1)_{n-k-2}l_k(x)=x(-A_kx+B_k)(x+1)\cdots(x+k-2)(x+1)\cdots(x+n-k-2)
$$$$
=x(-A_kx+B_k)\sum\limits_{j=0}^{q}e_{q-j}(\chi_k)x^j=B_ke_q(\chi_k)x+\sum\limits_{j=2}^{q+1}(B_ke_{q-j+1}(\chi_k)-A_ke_{q-j+2}(\chi_k))x^{j}
-A_kx^{q+2}
$$
$$
=B_ke_{p-1}(\chi_k)x+\sum\limits_{j=2}^{p}(B_ke_{p-j}(\chi_k)-A_ke_{p-j+1}(\chi_k))x^{j}
-A_kx^{p+1},
$$
where $2\leq{k}\leq{n/2}$, $q=n-4$, $p=n-3$,
$\chi_2=\{1,2,3,\ldots,n-4\}$ and
$$
\chi_k=\{1,1,2,2,\ldots,k-2,k-2,k-1,k,k+1,\ldots,n-k-2\},~~
k=3,4,\ldots...
$$
Note that each set $\chi_k$, $k=2,3,\ldots...$ has exactly $q=n-4$
elements.  If $k=1$ the formula is slightly different,
$$
\Phi_1(x)=B_1e_{p}(\chi_1)+\sum\limits_{j=1}^{p}(B_1e_{p-j}(\chi_1)-A_1e_{p-j+1}(\chi_1))x^{j}-A_1x^{p+1},~~\text{with}
$$$$
\chi_1=\{1,2,3,\ldots,n-3\}.
$$
The formula for $\Phi_k(x)$ shows that the coefficient at $x$ is
positive for all $k\geq{2}$ since $B_k>0$ for $0\leq{k}\leq{n/2}$
by its definition.  On the other hand, we know from
(\ref{eq:Phi0}) that the coefficient at $x$ is negative for $k=0$.
Hence, irrespective of the sign of the coefficient at $x$ for
$k=1$ our claim holds. Thus we can narrow our attention to the
coefficients at $x^j$ for $j=2,3,\ldots,n-2$.  Further, the
coefficients at $x^{n-2}$ are $-n(n-1),-A_1,
-A_2,\ldots,-A_{[n/2]}$. We have $A_k=A(k)$ for
$$
A(x)=n(n-1)-4x(n-x).
$$
Since $A(0)>0$, $A(n/2)<0$ and $A'(x)=8x-4n=0$ at $x=n/2$, $A(x)$
is decreasing on $[0,n/2]$ and changes sign exactly once.  So our
claim is true for the coefficients at $x^{n-2}$.

Finally we need to handle the general case of the coefficients at
$x^j$ for $j=2,3,\ldots,n-3$. It is easy to that
$\chi_{k-1}\prec^{W}\chi_k$ for $k=3,4,\ldots,[n/2]$ so that by
Lemma~1
$$
\frac{e_{p-j+1}(\chi_k)}{e_{p-j}(\chi_k)}<\frac{e_{p-j+1}(\chi_{k-1})}{e_{p-j}(\chi_{k-1})}
$$
for $j=2,3,\ldots,n-3$ and $k=3,4,\ldots,[n/2]$.  Further if
$A_k<0$ than it is clear that the coefficient at $x^{j}$ is
positive and there are no sign changes for such values of $k$.
Hence we take those $k$ for which $A_k\geq{0}$.  For such $k$ the
sequence $B_k/A_k$ is increasing, since
$$
\left(\frac{B(x)}{A(x)}\right)'=\frac{2n(n-1)(n-2x)}{A(x)^2}>0,~~B(x)=n(n-1)-2x(n-x).
$$
Now, if we assume that for some value of
$k\in\{3,4,\ldots,[n/2]\}$ the coefficient at $x^{j}$ in
$\Phi_k(x)$ is negative, i.e.
$$
B_ke_{p-j}(\chi_k)-A_ke_{p-j+1}(\chi_k)<0~\Leftrightarrow~\frac{B_k}{A_k}<\frac{e_{p-j+1}(\chi_k)}{e_{p-j}(\chi_k)}.
$$
Then for $k-1$ we will have
$$
\frac{B_{k-1}}{A_{k-1}}<\frac{B_k}{A_k}<\frac{e_{p-j+1}(\chi_k)}{e_{p-j}(\chi_k)}<\frac{e_{p-j+1}(\chi_{k-1})}{e_{p-j}(\chi_{k-1})}
~\Leftrightarrow~B_{k-1}e_{p-j}(\chi_{k-1})-A_{k-1}e_{p-j+1}(\chi_{k-1})<0,
$$
i.e. the coefficient at $x^j$ is again negative in
$\Phi_{k-1}(x)$.  This proves that there can be no more than one
change of sign in the sequence
$\{a_{2,j},a_{3,j},\ldots,a_{[n/2],j}\}$ for each
$j=2,3,\ldots,n-3$. It remains to consider $k=2$. Introduce
$$
\chi_2^{\varepsilon}=\{\varepsilon,1,2,\ldots,n-4\}.
$$
Clearly, $\chi_1\prec^W\chi_2^{\varepsilon}$ for each
$0<\varepsilon<1$ and $e_m(\chi_2^{\varepsilon})\to{e_m(\chi_2)}$
as $\varepsilon\to{0}$ for $m=0,1,\ldots$.  We have
$$
B_2e_{p-j}(\chi_2)-A_2e_{p-j+1}(\chi_2)<0~\Leftrightarrow~\frac{B_2}{A_2}<\frac{e_{p-j+1}(\chi_2)}{e_{p-j}(\chi_2)}~\Rightarrow~
\frac{B_2}{A_2}<\frac{e_{p-j+1}(\chi_2^{\varepsilon})}{e_{p-j}(\chi_2^{\varepsilon})}
$$
for sufficiently small $\varepsilon>0$ and
$$
\frac{B_1}{A_1}<\frac{B_2}{A_2}<\frac{e_{p-j+1}(\chi_2^{\varepsilon})}{e_{p-j}(\chi_2^{\varepsilon})}<\frac{e_{p-j+1}(\chi_1)}{e_{p-j}(\chi_1)}.~~\square
$$

\bigskip

\paragraph{4. Conjectures for higher order determinants.}  For $f(x;z)$ defined in (\ref{eq:f-defined}) let us
consider the Toeplitz determinant
$$
F_r(x,z)=\begin{vmatrix}
f(x;z) & f(x+1;z) & f(x+2;z) &\cdots & f(x+r-1;z)\\
f(x-1;z) & f(x;z) & f(x+1;z)&\cdots & f(x+r-2;z)\\
\vdots & \vdots & \vdots &\cdots &\vdots\\
f(x-r+1;z) & f(x-r+2;z) & f(x-r+3;z) &\cdots & f(x;z)\\
\end{vmatrix}.
$$
Compute
\begin{multline*}
F_r(x,z)=\begin{vmatrix}
\sum\limits_{k_1=0}^{\infty}f_{k_1}(x)_{k_1}\frac{z^{k_1}}{k_1!} & \sum\limits_{k_1=0}^{\infty}f_{k_1}(x+1)_{k_1}\frac{z^{k_1}}{k_1!} &\cdots & \sum\limits_{k_1=0}^{\infty}f_{k_1}(x+r-1)_{k_1}\frac{z^{k_1}}{k_1!}\\
\sum\limits_{k_2=0}^{\infty}f_{k_2}(x-1)_{k_2}\frac{z^{k_2}}{k_2!} & \sum\limits_{k_2=0}^{\infty}f_{k_2}(x)_{k_2}\frac{z^{k_2}}{k_2!} &\cdots & \sum\limits_{k_2=0}^{\infty}f_{k_2}(x+r-2)_{k_2}\frac{z^{k_2}}{k_2!}\\
\vdots & \vdots &\cdots &\vdots\\
\sum\limits_{k_r=0}^{\infty}f_{k_r}(x-r+1)_{k_r}\frac{z^{k_r}}{k_r!} & \sum\limits_{k_r=0}^{\infty}f_{k_r}(x-r+2)_{k_r}\frac{z^{k_r}}{k_r!} &\cdots & \sum\limits_{k_r=0}^{\infty}f_{k_r}(x)_{k_r}\frac{z^{k_r}}{k_r!}\\
\end{vmatrix}
\\[10pt]
=\sum\limits_{k_1,k_2,\ldots,k_r=0}^{\infty}f_{k_1}f_{k_2}\cdots{f_{k_r}}\frac{z^{k_1+k_2+\cdots+k_r}}{k_1!k_2!\cdots{k_r!}}
\begin{vmatrix}
(x)_{k_1} & (x+1)_{k_1} &\cdots & (x+r-1)_{k_1}\\
(x-1)_{k_2} & (x)_{k_2} &\cdots & (x+r-2)_{k_2}\\
\vdots & \vdots &\cdots &\vdots\\
(x-r+1)_{k_r} & (x-r+2)_{k_r} &\cdots & (x)_{k_r}\\
\end{vmatrix}
\\[10pt]
=\!\!\sum\limits_{n=0}^{\infty}\frac{z^n}{n!}\!\!\sum\limits_{k_1+k_2+\cdots+k_r=n}
\!\!\!\binom{n}{k_1,k_2,\ldots,k_r}f_{k_1}f_{k_2}\cdots{f_{k_r}}\!\begin{vmatrix}
(x)_{k_1} & (x+1)_{k_1} &\cdots & (x+r-1)_{k_1}\\
(x-1)_{k_2} & (x)_{k_2} &\cdots & (x+r-2)_{k_2}\\
\vdots & \vdots &\cdots &\vdots\\
(x-r+1)_{k_r} & (x-r+2)_{k_r} &\cdots & (x)_{k_r}\\
\end{vmatrix}.
\end{multline*}
Hence,
$$
F_r(x,z)=\sum\limits_{n=0}^{\infty}z^nP_n^r(x),
$$
where
$$
P_n^r(x):=\!\!\!\!\!\!\!\!\sum\limits_{k_1+k_2+\cdots+k_r=n}
\!\!\binom{n}{k_1,k_2,\ldots,k_r}f_{k_1}f_{k_2}\cdots{f_{k_r}}\begin{vmatrix}
(x)_{k_1} & (x+1)_{k_1} &\cdots & (x+r-1)_{k_1}\\
(x-1)_{k_2} & (x)_{k_2} &\cdots & (x+r-2)_{k_2}\\
\vdots & \vdots &\cdots &\vdots\\
(x-r+1)_{k_r} & (x-r+2)_{k_r} &\cdots & (x)_{k_r}\\
\end{vmatrix}.
$$
Of course, $P_n^2(x)=P_n(x)=Q_n^{1,1}(x-1)$. To conjecture a
reasonable  generalization of Theorem~\ref{th:size2} we need to
recall the notion of the P\'{o}lya frequency sequences, first
introduced by Fekete in 1912. They were studied in detail by
Karlin in \cite{Karlin}. The class of all P\'{o}lya frequency
sequences of order $1\leq{r}\leq\infty$ is denoted by $PF_r$ and
consists of the sequences $\{f_k\}_{k=0}^{\infty}$ such that all
minors of order $\leq{r}$ (all minors if $r=\infty$) of the
infinite matrix
\begin{equation}\label{eq:Polya}
\begin{bmatrix}
f_0 & f_1 & f_2 & f_3  &\cdots\\
0 & f_0 & f_1 & f_2 &\cdots \\
0 & 0 & f_0 & f_1 &\cdots \\
0 & 0 & 0 & f_0 &\cdots \\
\vdots & \vdots & \vdots &\vdots &\ddots \\
\end{bmatrix}
\end{equation}
are non-negative. Clearly,
$PF_1\supset{PF_2}\supset\cdots\supset{PF_{\infty}}$. The $PF_2$
sequences are precisely the log-concave sequences without internal
zeros. Our conjectures are
\begin{conjecture}\label{cj:big-pos}
Suppose $\{f_k\}_{k=0}^{n}\in{PF_{r}}$, $r\geq{2}$. Then the
polynomial $P_{n}^{r}(x)$ has degree $n-r(r-1)$ and positive
coefficients.
\end{conjecture}
\begin{conjecture}\label{cj:big-hurw}
Suppose $\{f_k\}_{k=0}^{n}\in{PF_{r}}$. Then the polynomial
$P_{n}^{r}(x)$ is Hurwitz stable.
\end{conjecture}
\begin{conjecture}\label{cj:big-zeros}
Suppose $\{f_k\}_{k=0}^{n}\in{PF_{\infty}}$. Then all zeros of the
polynomial $P_{n}^{r}(x)$ are real and negative for each
$r\geq{2}$.
\end{conjecture}
Again, Conjecture~\ref{cj:big-pos} follows from
Conjecture~\ref{cj:big-hurw} but both are independent of
Conjecture~\ref{cj:big-zeros}.

Conjectures~\ref{cj:zeros} and \ref{cj:big-zeros} bear certain
resemblance to the recent research of Br\"{a}nd\'{e}n
\cite{Braenden}, Grabarek \cite{Gr} and Yoshida \cite{Yoshida}.
Among other things, these works consider non-linear operators on
polynomials that preserve  the class of polynomials with real
negative zeros. According to the celebrated theorem of Aissen,
Schoenberg and Whitney \cite{ASW} the sequence  $\{f_0,
f_1,\ldots,f_n\}$ is a $PF_\infty$ sequence iff
$\sum_{k=0}^{n}f_kx^k$ has only real negative zeros.  In
particular, Br\"{a}nd\'{e}n found necessary and sufficient
conditions on the real sequence $\alpha_j$ to ensure that the
operators
\begin{equation}\label{eq:branden}
\{f_k\}_{k=0}^{n}\to\left\{\sum\limits_{j=0}^{\infty}\alpha_jf_{m-j}f_{m+j}\right\}_{m=0}^{n}~~\text{and}~~
\{f_k\}_{k=0}^{n}\to\left\{\sum\limits_{j=0}^{\infty}\alpha_jf_{m-j}f_{m+1+j}\right\}_{m=0}^{n-1}
\end{equation}
preserve $PF_{\infty}$. Here $f_i=0$ if
$i\notin\{0,1,\ldots,n\}$.  Using Br\"{a}nd\'{e}n's criterion
Grabarek showed in \cite{Gr} that the transformation ($p>0$ is an
integer)
\begin{equation}\label{eq:grabarek}
\{f_k\}_{k=0}^{n}\to\left\{\binom{2p-1}{p}f_m^2+\sum\limits_{j=1}^{p}(-1)^{j}\binom{2p}{p-j}f_{m-j}f_{m+j}\right\}_{m=0}^{n}
\end{equation}
preserves $PF_{\infty}$. Conjectures~\ref{cj:zeros} and
\ref{cj:big-zeros} also assert that certain non-linear
transformations preserve $PF_{\infty}$.  For $r=2$ this
transformation is easy to write explicitly.  Denote by $p_{n}(m)$,
$m=0,1,\ldots,n-2$, the coefficient at $x^m$ of the polynomial
$P_n(x)$.  Then
\begin{equation}\label{eq:p0}
p_n(0)=n!(f_1f_{n-1}-f_0f_{n}),
\end{equation}
and
$$
p_n(m)=\frac{1}{m!}\sum\limits_{0\leq{k}\leq{n/2}}f_kf_{n-k}\binom{n}{k}\frac{d^m}{dx^m}\Phi_k(x)_{|x=0},~~m=1,2,\ldots,n-2.
$$
For $k=0$ we have
$$
[\Phi_0(x)]^{(m)}_{|x=0}=-n(n-1)[(x+1)_{n-2}]^{(m)}_{|x=0}
=-n(n-1)\left[\sum\limits_{j=0}^{n-2}S^{n-1}_{j+1}x^j\right]_{|x=0}^{(m)}=-n(n-1)m!S^{n-1}_{m+1},
$$
where $S^p_j$ is unsigned Stirling numbers of the first kind that
can be defined by $(x)_p=\sum_{j=1}^{p}S^p_jx^j$.  For $k=1$ we
obtain,
$$
[\Phi_1(x)]^{(m)}_{|x=0}\!=\![(x+1)_{n-3}l_1(x)]^{(m)}_{|x=0}
\!=\!B_1[(x+1)_{n-3}]^{(m)}_{|x=0}-A_1m[(x+1)_{n-3}]^{(m-1)}_{|x=0}\!=\!m!(B_1S^{n-2}_{m+1}-A_1S^{n-2}_{m}),
$$
and for $2\leq{k}\leq{n/2}$, compute
$$
[\Phi_k(x)]^{(m)}_{|x=0}\!=\![(x)_{k-1}(x+1)_{n-k-2}l_k(x)]^{(m)}_{|x=0}
\!=\!B_k[(x)_{k-1}(x+1)_{n-k-2}]^{(m)}_{|x=0}-A_km[(x)_{k-1}(x+1)_{n-k-2}]^{(m-1)}_{|x=0},
$$
$$
[(x)_{k-1}(x+1)_{n-k-2}]^{(m)}_{|x=0}\!=\!\sum\limits_{i=0}^{m}\binom{m}{i}
\left[\sum\limits_{j=1}^{k-1}S^{k-1}_{j}x^j\right]_{|x=0}^{(i)}
\left[\sum\limits_{j=0}^{n-k-2}S^{n-k-1}_{j+1}x^{j}\right]_{|x=0}^{(m-i)}
$$
$$
=m!\sum\limits_{i=1}^{m}S^{k-1}_{i}S^{n-k-1}_{m-i+1},
~~~2\leq{k}\leq{n/2}.
$$
Hence,
$$
[\Phi_k(x)]^{(m)}_{|x=0}=B_km!\sum\limits_{i=1}^{m}S^{k-1}_{i}S^{n-k-1}_{m-i+1}
-A_km!\sum\limits_{i=1}^{m-1}S^{k-1}_{i}S^{n-k-1}_{m-i},~~~2\leq{k}\leq{n/2}.
$$
Finally, we get for $m=1,2,\ldots,n-2$,
\begin{multline}\label{eq:pm}
p_n(m)=-n(n-1)S^{n-1}_{m+1}f_0f_{n}+f_1f_{n-1}n(B_1S^{n-2}_{m+1}-A_1S^{n-2}_{m})
\\
+\sum\limits_{2\leq{k}\leq{n/2}}f_kf_{n-k}\binom{n}{k}\left(B_k\sum\limits_{i=1}^{m}S^{k-1}_{i}S^{n-k-1}_{m-i+1}
-A_k\sum\limits_{i=1}^{m-1}S^{k-1}_{i}S^{n-k-1}_{m-i}\right)
\\
=-n(n-1)S^{n-1}_{m+1}f_0f_{n}+\sum\limits_{1\leq{k}\leq{n/2}}f_kf_{n-k}\binom{n}{k}
\sum\limits_{i=0}^{m}S^{k-1}_{i}\left(B_kS^{n-k-1}_{m-i+1}-A_kS^{n-k-1}_{m-i}\right),
\end{multline}
where
$$
S^{p}_{q}=0,~~q>p,~~S_0^{p}=0,~~p\geq{1},~~S_{0}^{0}=1.
$$
So Conjecture~\ref{cj:zeros} can be restated as the assertion that
the non-linear operator
$$
\{f_k\}_{k=0}^{n}\to
\left\{\sum\limits_{0\leq{j}\leq{n/2}}f_jf_{n-j}P_{j,m}\right\}_{m=0}^{n-2},
$$
where the numbers $P_{j,m}$ can be read off (\ref{eq:p0}) and
(\ref{eq:pm}), preserves $PF_{\infty}$.  Both Br\"{a}nd\'{e}n's
transformation (\ref{eq:branden}) our transformation above are
bilinear forms but of somewhat different character. One may ask
then what conditions on the numbers $P_{k,m}$ would ensure the
preservation of $PF_{\infty}$.

\paragraph{5. Some remarks on numerical experiments.}  In order to
run numerical experiments with Conjectures~1 to 6 one has to be
able to generate $PF_r$ sequences.  For $r=2$ and $r=\infty$ the
methods are quite clear. Setting
$$
\delta_k=\frac{f_k^2}{f_{k-1}f_{k+1}}
$$
we obtain for $\{f_k\}_{k=0}^{\infty}\in{PF_2}$:
\begin{equation}\label{eq:fdelta}
f_k=f_0^{k+1}\delta_1^k\delta_2^{k-1}\cdots\delta_k,
\end{equation}
where $f_0>0$ and $0<\delta_j\leq{1}$, $j=1,2,\dots,n$,
$\delta_j=0$, $j>n$. Hence, we can parameterize all $PF_2$
sequences by sequences with elements from $(0,1]$. Generating the
latter randomly we get a random $PF_2$ sequence. Next, for
$r=\infty$ we can simply generate $n$ random positive numbers
$a_1,a_2,\ldots,a_n$ and compute the coefficients of the
polynomial $\prod_{i=1}^{n}(x+a_i)$ producing by
Aissen-Shoenberg-Whitney theorem a $PF_{\infty}$ sequence.
According to the same theorem all finite $PF_{\infty}$ sequences
are obtained in this way.

The situation is less clear for $3\leq{r}<\infty$.  I am unaware
of any method to parameterize all $PF_r$ sequences for these
values of $r$.  However, some subclasses can be parameterized. One
possible method is provided by the following result of Katvova and
Vishnyakova \cite[Corollary of Theorem~5]{KV}: if nonnegative
sequence $\{f_n\}_{n=0}^{\infty}$ satisfies
$$
f_n^2\geq
4\cos^2\left(\frac{\pi}{r+1}\right)f_{n-1}f_{n+1},~~~n\geq{1},
$$
then $\{f_n\}_{n=0}^{\infty}\in{PF_{r}}$. This implies that if we
choose $0<\delta_j\leq\left(4\cos^2\frac{\pi}{r+1}\right)^{-1}$
then the sequence generated by (\ref{eq:fdelta}) is a $PF_{r}$
sequence.  Another method to produce a finite $PF_r$ sequence
follows from Shoenberg's theorem \cite{Sch} stating that the
coefficients of a polynomial with zeros lying in the sector
$|\arg{z}-\pi|<\pi/(r+1)$  form a $PF_r$ sequence. Hence,
generating such zeros randomly and doubling their number by adding
the complex conjugate to each we get a polynomial with  $PF_{r}$
coefficients.

Finally, Ostrovskii and Zheltukhina \cite{OsZh} parameterized a
large subclass of $PF_3$ sequences.  Namely, a $PF_3$ sequence
$\{f_0,f_1,f_2,\ldots,\}$ is $Q_3$ if all truncated sequences
$\{f_i\}_{i=0}^{n}$ are also $PF_3$ for each $n=1,2,\ldots$. The
main Theorem of \cite{OsZh} states that a sequence
$\{f_0,f_1,f_2,\ldots,\}$ is $Q_3$ iff  $f_0>0$,
$f_1=f_0\beta\geq{0}$ and
$$
f_n=\frac{f_0\beta^n\delta_2^{n-1}\delta_3^{n-2}\cdots\delta_{n-1}^2\delta_n}{\alpha_2^{n/2}\alpha_3^{(n-1)/2}\alpha_4^{(n-2)/2}\cdots\alpha_{n-1}^{3/2}\alpha_n},
$$
where
$$
\alpha_2=1+\delta_2,~\alpha_3=1+\delta_3\sqrt{\alpha_2},~\alpha_4=1+\delta_4\sqrt{\alpha_3},\ldots,
~~0\leq\delta_j\leq{1},~~j=2,3,\ldots
$$
and the sequence $\{\delta_j\}$ has no internal zeros.  This
theorem provides a simple method of generating random $Q_3$
sequences.

\paragraph{6. Acknowledgements.} I thank Lukasz Grabarek, Sergei Kalmykov, Mikhail
Tyaglov and Dennis Stanton for numerous useful discussions
concerning conjectures~\ref{cj:positivity} and \ref{cj:hurwitz}.

\end{document}